# Properties of Subsets of Unobserved Items Constructed by Using Blackwell's Bet

By James D. Stein
Department of Mathematics and Statistics, California State University (Long Beach)

**Abstract –** Blackwell's Bet (also known as the Cover Guessing Game) is an extremely clever use of a random variable to improve the probability of guessing which of two numbers is larger when you are only allowed to observe one of those numbers. This technique has hitherto been used to guess or predict outcomes with better-than-chance results. We show that this technique can also be used to construct subsets of unobserved items which have expected values below and above the expected value of the parent set. We discuss this in two situations – when the items in question can have arbitrary values, and when the items in question are Bernoulli trials. We also include simulations in an instance when the actual formulas require extensive.iteration of recursively-defined relations.

**Blackwell's Bet**

Suppose that two envelopes contain differing amounts of money $S < L$. You are allowed to open one envelope and count the money in it. You are then offered the opportunity to either keep the money in that envelope or take the money in the unopened envelope. Blackwell ([1]) showed that you can improve your probability of successfully choosing the larger amount with the following strategy.

Assume that $S < L$, and we open an envelope containing D dollars, where D is equally likely to be S or L. Let r be a random variable. If $r < D$, keep the money in the opened envelope; if not, take the money in the unopened envelope. Using this strategy, the probability of a successful guess is ½ + ½ the probability thar r is located between S and L. The proof is simple and can be found on p.3 of [4].

A variant of this technique can also be used to guess which of two Bernoulli trials has a higher success probability with just a single observation. Choose one of the two trials at random, and conduct a single observation of it. If the observation is a success, guess that you are looking at the trial with the higher success probability, and conversely. The probability of a successful guess in this case is ½ + ½ the difference between the larger and smaller success probabilities of the Bernoulli trials.

The Blackwell technique involves selecting two items and choosing one equiprobably to examine. This process tells us something about the unexamined item – and it is that information which is the subject of this paper.

**Section I – Creating Subsets of Unexamined Items with Different Expected Values**

If one were faced with the problem of guessing which of two marbles is the heavier, one can implement the strategy in Blackwell's Bet without knowing either the weight of the chosen marble or the value of a random variable. Simply take a bag of marbles of assorted weights and use a balance scale to compare the weights of the chosen marble and a random marble selected from the bag. If the random marble weighs less than the chosen marble, guess that the chosen marble is the heavier, and vice-versa. The mathematics is the same, and the probability of guessing which marble is the heavier is ½ + one-half the probability that the weight of the randomly chosen marble lies between the weights of the two given marbles. As a result, using Blackwell's Bet might result in obtaining information at less expense than if precise measurements, such as weighing a marble, were used. Unless otherwise noted throughout the remainder of this paper, we will discuss the problems involved by using the example of sorting marbles in a bag by weight. This is done for the purpose of clarity – using expressions such as "weight of a marble" instead of "value of a random variable" makes things more understandable. Or so it appears to the author

Additionally, Blackwell's Bet works no matter what random variable is used, as long as the probability that the random variable lies between the lower and higher numbers is non-zero. Tailoring random variables to specific situations often results in better guessing probabilities or other desirable outcomes, as will be shown in this paper.

**Creating Subsets with Expected Values Different from the Parent Set**

Suppose you have a bag containing marbles of varying weights and a scale. The goal is to create two subsets A and B of marbles, such that the expected value of the weight of a marble in A is less than the expected value of the weight of a marble from the entire bag, and the expected value of the weight of a marble from the entire bag is less than the expected value of the weight of a marble in B, **without examining in any way the marbles that are assigned to A or B**. You can weigh or compare other marbles, but you cannot weigh or compare any marble you assign to A or B.

It is possible to do this so that in many situations both the expected weight of a marble in A and the expected weight of a marble in B satisfy the desired conditions.

**Sorting Unexamined Items**

Blackwell's Bet suggests the following strategy – take pairs of marbles from the bag, examine one of those two marbles equiprobably, and use a random variable as in Blackwell's Bet to guess which is larger. If the Blackwell's Bet technique results in guessing that the unexamined one is the smaller of the two items, assign it to A, otherwise assign it to B.

The following Lemma gives the result of this process, which we will call the standard assignment procedure (SAP), in the case of two marbles with different weights.

**Lemma 1** – Suppose that we are confronted with two marbles with unknown weights $S < L$. We choose one of the two marbles equiprobably and measure (or compare) it with a randomly selected weight. Call this weight X (the weight of the eXamined marble) and the weight of the other marble U (the weight of the Unexamined marble). Now choose a random variable r from a distribution with $P(x<S) = p$ and $P(x>L) = q$. We write $U \rightarrow A$ to indicate that U is assigned to A. The decision rule is that if $r < X$, $U \rightarrow A$, and if $r > X$, $U \rightarrow B$.

Then $E(U\rightarrow A) = (Lp + S(1-q))/(1+p-q) < (L+S)/2 < E(U\rightarrow B) = (L(1-p) + qS)/(1-p+q)$, where $E(U\rightarrow A)$ is the expected value of the weight of an unexamined marble assigned to A and $E(U\rightarrow B)$ is defined similarly.

**Proof:** We envision a bag containing equal quantities of marbles with weights S and L – although we do not know what the numerical values of these weights are. The marbles are extracted from the bag in pairs, one of each weight. Note that X and U are numerical quantities, one of which will be L and one of which will be S – but comparison allows one to assign U to either A or B without knowing the numerical value of either X or r.

Let $P(U\rightarrow A)$ denote the probability that $U \rightarrow A$. We can think of $P(U\rightarrow A)$ as the fraction of a marble that is assigned to A. Let $W(U\rightarrow A)$ denote the total weight of the fractional marble assigned to A. Similar definitions are made for $U \rightarrow B$.

With probability ½, $X=S$ (so $U=L$). In this case, $P(U \to A) = P(X \to B) = P(S \to B) = p$. So the combined probability of $X = S$ and $U \to A = 0.5p$. The combined probability of $X = S$ and $U \to B$ is $0.5(1-p)$. So $W(U \to A) = 0.5Lp$ and $W(U \to B) = 0.5L(1-p)$.

With probability ½, $X=L$ (so $U=S$). In this case, $P(U \to A) = P(X \to B) = P(L \to B) = 1-q$. So the combined probability of $X = L$ and $U \to A = 0.5(1-q)$. The combined probability of $X = L$ and $U \to B$ is $0.5q$. So $W(U \to A) = 0.5S(1-q)$ and $W(U \to B) = 0.5Sq$.

So $P(U \to A) = 0.5(p+1-q)$ and $W(U \to A) = 0.5(Lp+S(1-q))$. Therefore, thinking of the probability $P(U \to A)$ as the fraction of the unexamined marble assigned to A, we see that $E(U \to A) = (Lp + S(1-q))/(1+p-q)$, and $P(U \to B) = 0.5((1-p)+q)$ and $W(U \to B) = 0.5(L(1-p)+Sq)$. Similarly, $E(U \to B) = (L(1-p) + Sq)/(1-p+q)$.

Notice that the mean value of L and S is $(L+S)/2$. We now show the inequality $E(U \to A) < (L+S)/2$

In order to show $E(U \to A) < (L+S)/2$, by cross-multiplying this is true if

$$2(Lp+S(1-q)) < (1+p-q)(L+S)$$

Expanding and simplifying

$$2Lp + 2S - 2Sq < L + S + Lp - Lq + Sp - Sq$$

$$S(1-p-q) < L(1-p-q)$$

The latter equality is true whenever $p+q < 1$ and is an equality if $p+q=1$. Of course, it is impossible that $p+q > 1$.

The proof of $(L+S)/2 < E(U \to B)$ is similar. Proving it requires that we show

$(L+S)/2 < (L(1-p)+qS)/(1-p+q)$. By cross-multiplying this is true if

$$(L+S)((1-p+q) < 2(L(1-p) +qS)$$

Expanding and subtracting $L(1-p) + qS$ from both sides shows that this is equivalent to

$$Lq + S(1-p) < L(1-p) + Sq$$

$$S(1-p-q) < L(1-p-q), \text{ as before. } \blacksquare$$

**Cor. 1.1 –** Assume that $0 < S < L < 1$ and that the uniform probability distribution on [0,1] is used to select the random variable r. Then $p = S$, $q = 1-L$, and

$P(U\to A) = 0.5(S+L)$ $W(U\to A) = 0.5(2LS)$ $E(U\to A) = 2LS/(S+L)$

$P(U\to B) = 0.5(2-S-L)$ $W(U\to B) = 0.5((1-L)S+(1-S)L)$ $E(U\to B) = (L+S-2LS)/(2-S-L)$

Note that if AM(x,y) denotes the arithmetic mean of x and y, and GM(x,y) denotes the geometric mean of x and y, then

$$E(U\to A) = GM(S,L)^2/AM(S,L) \text{ and } E(U\to B) = (AM(S,L)-GM(S,L)^2)/(1-AM(S,L))$$

These expressions hold only when $0 < S < L < 1$ and we use the uniform Blackwell random variable. To see this, consider the problem of solving the two equations (in the unknowns p and q)

(1) $(Lp + S(1-q))/(1+p-q) = 2SL/(S+L)$

(2) $(L(1-p) + Sq)/(1-p+q) = (L+S-2SL)/(2-S-L)$

Equation (1) simplifies to $Lp + Sq = S$, and equation (2) simplifies to $(L-1)p + (S-1)q = L-1$. The unique solution to these equations is $p = S$, $q = 1 – L$, and since $q \geq 0$, $L \leq 1$.

However, if we let p = q = S/(L+S), then we obtain

$$E(U \to A) = 2LS/(L+S) = GM(S,L)^2/AM(S,L)$$

$$\text{and } E(U \to B) = (L^2 + S^2)/(L+S) = 2AM(S,L) - GM(S,L)^2)/(AM(S,L))$$

**Cor. 1.2 –** Assume that 1 < S < L, but we use ln S and ln L as the weights of the marbles, and that we choose a random variable r from a distribution with P(x< ln S) = p and P(x> ln L) = q. Then

$$E(U \to A) = \frac{p \ln L + (1-q) \ln S}{1+p-q} = \ln(L^p S^{1-q})^{\frac{1}{1+p-q}} < \frac{\ln S + \ln L}{2} = \ln(\sqrt{LS}) \text{ and}$$

$$\ln(\sqrt{LS}) < \frac{(1-p)\ln L + q \ \ln S}{1+q-p} = \ln(L^{1-p} S^q)^{\frac{1}{1+q-p}} = \ E(U \to B)$$

**Section II –When All Marbles Have Different Weights in (0,1)**

**W**e assume that the marbles have different weights $W_1 < W_2 < \ldots < W_n$ between 0 and 1. We use SAP and the uniform distribution on (0,1) to select the random value for comparison. In the case where we only have one marble of each weight, there are two separate cases – where the bag contains even and odd numbers of marbles. In this case, the sorting procedure for the bag containing an odd number of marbles leaves a single unexamined marble in the bag at the end – this marble is assigned to A or B randomly with probability ½. Doing this will not affect the conclusion of the following Theorem 1.

The case where each pair consists of marbles of differing weights includes such important situations as when we are confronted with ordinal data, and has the advantage of being relatively easy to analyze. As in Section I, X is the weight of the examined marble and U the weight of the unexamined marble.

**Theorem 2.1** – Let m be the expected value of the weight of the marbles. Then $E(U \to A) < m < E(U \to B)$.

**Proof:** Let P be the Blackwell random variable, and define $p_k = P([0, W_k))$ for k = 1, 2 ,.., n. Since the weights of the marbles are assumed to increase, the values of $p_k$ are non-decreasing.

Assume we have chosen the pair $\{W_i, W_j\}$ for examination. With probability ½, $X = W_i$. The combined probability that $X = W_i$ and $X \to B$ is ½ $p_i$; this is also the combined probability that $U = W_j \to A$. In this case, the expected value of the weight of the marble assigned to A is ½ $p_i W_j$. Similarly, the combined probability that $X = W_j$ and $X \to B$ is ½ $p_j$; this is also the combined probability that $U = W_i \to A$ and the expected value of the weight of the marble assigned to A is ½ $p_j W_i$. Therefore, the probability of assigning the unexamined marble to A is ½$(p_i + p_j)$ and the expected value of the weight of the unexamined marble assigned to A is ½$(p_i W_j + p_j W_i)$.

Since each pair of marbles occurs with probability 2/(n(n-1)), we have

$$P(U \to A) = \frac{1}{n(n-1)} \sum_{i=1}^{n-1} \sum_{j=i+1}^{n} (p_i + p_j)$$

$$W(U \to A) = \frac{1}{n(n-1)} \sum_{i=1}^{n-1} \sum_{j=i+1}^{n} (p_i W_j + p_j W_i)$$

$$E(U \to A) = \frac{\sum_{i=1}^{n-1} \sum_{j=i+1}^{n} (p_i W_j + p_j W_i)}{\sum_{i=1}^{n-1} \sum_{j=i+1}^{n} (p_i + p_j)}$$

The expected value of the weight of the marbles in the bag is $m = (W_1 + \ldots + W_n)/n$. To show that $E(U \to A) < m$, it therefore suffices to show that

$$n \sum_{i=1}^{n-1} \sum_{j=i+1}^{n} (p_i W_j + p_j W_i) < (W_1 + \ldots + W_n) \sum_{i=1}^{n-1} \sum_{j=i+1}^{n} (p_i + p_j)$$

The double sum in the right side of the above inequality simplifies to

$(n-1)\sum_{i=1}^{n} p_i$. We therefore must show that

$n \sum_{i=1}^{n-1} \sum_{j=i+1}^{n} (p_i W_j + p_j W_i)$ <(n-1) ($W_1$ + … + $W_n$) $\sum_{i=1}^{n} p_i$

Each of the terms of the form $p_iW_j$ with i≠j appears n times on the left and n-1 times on the right. Subtracting (n-1) times each term of the form $p_iW_j$ with i≠j from both sides reduces the problem to showing that

$$\sum_{i=1}^{n-1} \sum_{j=i+1}^{n} (p_i W_j + p_j W_i) < (n-1) \sum_{i=1}^{n} p_i W_i$$

The expression on the left is simply the sum of all $p_iW_j$ with i≠j.

Observe that, for j > i, $(p_j - p_i)(W_j - W_i) > 0$. Therefore

$$0 < \sum_{i=1}^{n-1} \sum_{j=i+1}^{n} (p_i W_i - p_j W_i - p_i W_j + p_j W_j)$$

Adding all the terms of the form $p_iW_j$ with i≠j to both sides yields the desired result.

Showing that m < E(U→B) is the same proof by switching A and B, replacing $p_i$ by 1- $p_i$ and reversing the direction of the inequality. █

We now tackle a more general situation – when we allow pairs of marbles with the same weight.

**Section III – Allowing Pairs of Marbles with Equal Weights**

If we allow differing probabilities for the pairs {S,S}, {S,L}, and {L,L}, the SAP results in situations in which A and B have different expected values from the overall expected value.

**Theorem 3.1** - Suppose that the pairs {S,S}, {L,L}, and {S,L} occur with probabilities α, β and 1 - α – β, and suppose that P(x<S) = p and P(x>L) = q. Then the overall expected value m is given by

$$m = \alpha S + \beta L + (1-\alpha-\beta)\frac{S+L}{2}$$

The expected value of the weights of the unobserved marbles assigned to A and B are given by

$$E(U \to A) = \frac{\alpha p S + \beta(1-q)L + (1-\alpha-\beta)\frac{(Lp + S(1-q))}{2}}{\alpha p + \beta(1-q) + (1-\alpha-\beta)\frac{p+1-q}{2}}$$

$$E(U \to B) = \frac{\alpha(1-p)S + \beta q L + (1-\alpha-\beta)\frac{(L(1-p) + Sq)}{2}}{\alpha(1-p) + \beta q + (1-\alpha-\beta)\frac{1-p+q}{2}}$$

**Proof:** We construct the following table.

| **Pair** | **Probability of Appearance** | **P(U→A)** | **W(U→A)** |
|---|---|---|---|
| {S,S} | α | p | pS |

With probability 1, S will be examined. It will be assigned to B with probability p, which is also the probability that U=S will be assigned to A.

| **Pair** | **Probability of Appearance** | **P(U→A)** | **W(U→A)** |
|---|---|---|---|
| {L,L} | β | 1-q | (1-q)L |

With probability 1, L will be examined. It will be assigned to B with probability 1-q, which is also the probability that U=L will be assigned to A.

| {S,L} | 1-α-β | (p+1-q)/2 | (Lp+S(1-q))/2 |
|---|---|---|---|

The last two numbers are the conclusion of Lemma 1.

If we randomly select a pair of marbles, and then randomly choose one from the selected pair, the expected value of the weight of the selected marble from the pair {S,S} is S, from the pair {L,L} is L, and from the pair {S,L} is (S+L)/2. So the expected value of the weight of a marble in the bag is

$$m = \alpha S + \beta L + (1 - \alpha - \beta)\frac{S + L}{2}$$

From the above table, the expected value of weight of the unexamined marbles assigned to A is given by

$$E(U \to A) = \frac{\alpha p S + \beta(1 - q)L + (1 - \alpha - \beta)\frac{(Lp + S(1 - q))}{2}}{\alpha p + \beta(1 - q) + (1 - \alpha - \beta)\frac{p + 1 - q}{2}}$$

We can similarly construct the following table.

| Pair | Probability of Appearance | P(U→B) | W(U→B) |
|---|---|---|---|
| {S,S} | α | 1-p | (1-p)S |

With probability 1, S will be examined. It will be assigned to A with probability 1-p, which is also the probability that U=S will be assigned to B.

| Pair | Probability of Appearance | P(U→B) | W(U→B) |
|---|---|---|---|
| {L,L} | β | q | qL |

With probability 1, L will be examined. It will be assigned to A with probability q, which is also the probability that U=L will be assigned to B.

| Pair | Probability of Appearance | P(U→B) | W(U→B) |
|---|---|---|---|
| {S,L} | 1-α-β | (1-p+q)/2 | (L(1-p)+qS))/2 |

See proof of Lemma 1.

The overall expected value is given by

$$m = \alpha S + \beta L + (1 - \alpha - \beta)\frac{S + L}{2}$$

As above, the expected value of the weight of the unexamined marbles assigned to B is given by

$$E(U \rightarrow B) = \frac{\alpha(1 - p)S + \beta qL + (1 - \alpha - \beta)\frac{(L(1 - p) + Sq)}{2}}{\alpha(1 - p) + \beta q + (1 - \alpha - \beta)\frac{1 - p + q}{2}} \quad \blacksquare$$

We now examine some interesting special cases in which we can easily determine how the expected values $E(U \rightarrow A)$ and $E(U \rightarrow B)$ compare with the overall expected value m.

**Cor. 3.1 - A**ssume that $\alpha = \beta$. Then $E(U \rightarrow A) < m$ and $E(U \rightarrow B) > m$ iff $\alpha < ¼$.

**Proof:** From Theorem 3.1, we obtain by substitution

$$m = \alpha(S + L) + (1 - 2\alpha)\frac{S + L}{2} = \frac{S + L}{2}$$

numerator of $E(U \rightarrow A) = \alpha(pS + (1 - q)L) + (1 - 2\alpha)\frac{(Lp + S(1 - q))}{2}$

denominator of $E(U \rightarrow A) = \alpha\big(p + (1 - q)\big) + (1 - 2\alpha)\frac{p + 1 - q}{2} = \frac{p + 1 - q}{2}$

We wish to determine when $E(U \rightarrow A) < m$. To do so, take the above expressions for the numerator and denominator and get a simplified expression for $E(U \rightarrow A)$ by multiplying both the numerator and denominator by 2, and then dividing the numerator by the denominator. Take this expression, assume $E(U \rightarrow A) < m$, and multiply both sides by 2(1+p-q). This results in the following inequality.

$$4\alpha(1 - q - p)(L - S) + 2Lp + 2S(1 - q) < (1 + p - q)(S + L)$$

Subtract 2Lp + 2S(1-q) from both sides. The above inequality simplifies to

$$4\alpha(1 - q - p)(L - S) < (1 - p - q)(L - S)$$

Since L > S and 1 > p + q, the above simplifies to 4α < 1, with equality holding when α = ¼. This is the case when choosing a marble S has the same probability of choosing a marble L, and both choices are independent.

We now wish to determine when $E(U \rightarrow B) > m$. To do so, take the simplified expressions for the numerator and denominator and get a simplified expression for $E(U \rightarrow B)$ by multiplying both the numerator and denominator by 2, and then dividing the numerator by the denominator. Take this expression, assume $E(U \rightarrow B) > m$, and multiply both sides by 2(1-p+q). This results in the following inequality.

$$4\alpha(1 - p - q)(S - L) + 2L(1 - p) + 2Sq > (1 - p + q)(S + L)$$

Subtract 2L(1-p) + 2Sq from both sides. The above inequality simplifies to

$$4\alpha(1 - p - q)(S - L) > 1 - p - q)(S - L)$$

Since L > S and 1 > p + q, the above simplifies to 4α < 1, with equality holding when α = ¼. As before, this is the case when choosing a marble S has the same probability of choosing a marble L, and both choices are independent. █

Cor. 3.1 raises the possibility that $E(U \to A) = m$ if the selection process for the two marbles results in independent choices of S and L. We show this is the case.

**Cor. 3.2** – Suppose that the choice of marbles with weights S and L is independent. Then $E(U \to A) = m$.

**Proof:** If the choice of marbles is independent, then $\alpha = a^2$, $\beta = b^2$ where $a + b = 1$. We have $1 - \alpha - \beta = 1 - a^2 - b^2 = (1-a)(1+a) - b^2 = b(1+a) - b^2 = b(1+a-b) = 2ab$. We also have $m = a^2S + b^2L + ab(S+L) = a(a+b)S + b(a+b)L = aS + bL$.

We now proceed to simplify the numerator and denominator of $E(U \to A)$.

$$\text{numerator of } E(U \to A) = a^2pS + b^2(1-q)L + ab(Lp + S(1-q))$$

$$= S(a^2p + ab(1-q)) + L(b^2(1-q) + abp)$$

$$= aS(ap + b(1-q)) + bL(b(1-q) + ap)$$

$$= (ap + b(1-q))(aS + bL) = (ap + b(1-q))m$$

$$\text{denominator of } E(U \to A) = a^2p + b^2(1-q) + ab(p + (1-q))$$

$$= ap(a + b) + b(1-q)(b + a)$$

$$= ap + b(1-q)$$

So $E(U \to A) = m$. ■

It is also true that if the choice of marbles with weights S and L is independent, $E(U \to B) = m$. The proof is similar and will be omitted.

We can see that, in general, this is the only case for which $E(U \to A) = m$ by taking as the Blackwell random variable the uniform distribution on [0,1], where $0 = S < L < 1$. In this case, p $= S$ and $1 - q = L$. Then $E(U \rightarrow A) = \frac{\beta L^2}{\frac{1-\alpha+\beta}{2}L} = \frac{2\beta}{1-\alpha+\beta}L$ and $m = \frac{1-\alpha+\beta}{2}L$. Equating these two expressions, we see that $4\beta = (1 - \alpha + \beta)^2$. The only solution to this is $\beta = (1 - \sqrt{\alpha})^2$, which is seen to be the case having $\alpha = a^2$, $\beta = b^2$ where $a + b = 1$.

**Ensuring E(U→A) < m < E(U→B) Given Independent Selection of Marbles with Weights S and L**

Although Cor. 3.2 says this cannot be done using the SAP for two marbles, Cor. 3.1 suggests a way around this – by making sure that we reduce the relative frequency of examining the pairs {L,L} and {S,S}.

When we are examining marbles, we can place each marble on one side of a balance scale without observing the marbles. One way to do this is simply to drop each marble into a container with opaque sides. If the balance scale shows that the two are equal, examine neither. If they are unequal, remove both marbles and drop them simultaneously into a container with opaque sides, and then remove and observe one of the marbles using the SAP. Similar procedures can be adopted for differing physical parameters.

When given a pair {x,y} of numerical values, if it is possible to compute $(x-y)^2$ or $|x-y|$ without observing x or y (which could be done blind by passing the numbers to a computer program which performs this function), one could discard the pair {x,y} without observing either if it is shown that $(x-y)^2$ or $|x-y| = 0$.

The cost of doing this in either the physical or numerical case is that fewer items will be assigned to A and B.

It might be possible to employ other versions of Blackwell's Bet to ensure that the subsets A and B have differing expected values. The author tried two schemes without success. The first was to use a variant of King of the Hill, selecting three marbles and observing two of them. The second was a version of a tournament, in which four marbles were selected and competed in a standard tournament format, in which three of the four were observed. Neither was successful in achieving the desired result, leaving the author to speculate whether the creation of the subsets A and B were possible if the original set contained infinitely many marbles.

**Section IV - Recursion Relations for Finite Collections of Two Bernoulli Trials**

After $k = j_1 + j_2 + j_3 + j_4$ unobserved marbles have been assigned, let $P(j_1,j_2,j_3,j_4)$ denote the probability that $j_1$ marbles of type S have been assigned to A, $j_2$ marbles of type S have been assigned to B, $j_3$ marbles of type L have been assigned to A, and $j_4$ marbles of type L have been assigned to B.

(1) To compute $P(j_1+1,j_2,j_3,j_4)$, either the pair {S,S} or {S,L} must have been drawn.

If {S,S} has been drawn, then $j_1 + j_2 < N_S - 1$. If $j_1 + j_2 \geq N_S - 1$, let $\alpha = 0$. If $j_1 + j_2 < N_S - 1$, then the probability of drawing {S,S} is $C(N_S - j_1 - j_2,2)/C(N_S + N_L - k,2)$. The probability that S is observed and the unobserved S is assigned to A is S. Then let

$\alpha = S\ C(N_S - j_1 - j_2,2)/C(N_S + N_L - k,2)$.

If {S,L} has been drawn, then $j_1 + j_2 < N_S$ and $j_3 + j_4 < N_L$. If either $j_1 + j_2 = N_S$ or $j_3 + j_4 = N_L$ let $\beta = 0$. If not, the probability of {S,L} is $(N_S - j_1 - j_2)(N_L - j_3 - j_4)/(N_S + N_L - k,2) = \gamma$. The probability that L is observed and S is assigned to A is 0.5L. Then let $\beta = 0.5L\gamma$.

Then $P(j_1+1,j_2,j_3,j_4) = (\alpha + \beta)\ P(j_1,j_2,j_3,j_4)$.

(2) To compute $P(j_1,j_2+1,j_3,j_4)$, either the pair {S,S} or {S,L} must have been drawn.

If {S,S} has been drawn, then $j_1 + j_2 < N_S - 1$. If $j_1 + j_2 \geq N_S - 1$, let $\alpha = 0$. If $j_1 + j_2 < N_S - 1$, then the probability of drawing {S,S} is $C(N_S - j_1 - j_2,2)/C(N_S + N_L - k,2)$. The probability that S is observed and the unobserved S is assigned to B is 1-S. Then let

$\alpha = (1\text{-}S)\ C(N_S - j_1 - j_2,2)/C(N_S + N_L - k,2)$.

If {S,L} has been drawn, then $j_1 + j_2 < N_S$ and $j_3 + j_4 < N_L$. If either $j_1 + j_2 = N_S$ or $j_3 + j_4 = N_L$ let $\beta = 0$. If not, the probability of {S,L} is $(N_S - j_1 - j_2)(N_L - j_3 - j_4)/(N_S + N_L - k,2) = \gamma$. The probability that L is observed and S is assigned to B is 0.5(1-L). Then let $\beta = 0.5(1\text{-}L)\gamma$.

Then $P(j_1,j_2+1,j_3,j_4) = (\alpha + \beta)\ P(j_1,j_2,j_3,j_4)$.

(3) To compute $P(j_1,j_2,j_3+1,j_4)$, either the pair {L,L} or {S,L} must have been drawn.

If {L,L} has been drawn, then $j_3 + j_4 < N_L - 1$. If $j_3 + j_4 \geq N_L - 1$, let $\alpha = 0$. If $j_3 + j_4 < N_L - 1$, then the probability of drawing {L,L} is $C(N_L - j_3 - j_4,2)/C(N_S + N_L - k,2)$. The probability that L is observed and the unobserved L is assigned to A is L. Then let

$\alpha = L\ C(N_L - j_3 - j_4,2)/C(N_S + N_L - k,2)$.

If {S,L} has been drawn, then $j_1 + j_2 < N_S$ and $j_3 + j_4 < N_L$. If either $j_1 + j_2 = N_S$ or $j_3 + j_4 = N_L$ let $\beta = 0$. If not, the probability of {S,L} is $(N_S - j_1 - j_2)(N_L - j_3 - j_4)/(N_S + N_L - k,2) = \gamma$. The probability that S is observed and L is assigned to A is 0.5S. Then let $\beta = 0.5S\gamma$.

Then $P(j_1,j_2,j_3+1,j_4) = (\alpha + \beta)\ P(j_1,j_2,j_3,j_4)$.

(4) To compute $P(j_1,j_2,j_3,j_4+1)$, either the pair {L,L} or {S,L} must have been drawn. If {L,L} has been drawn, then $j_3 + j_4 < N_L - 1$. If $j_3 + j_4 \geq N_L - 1$, let $\alpha = 0$. If $j_3 + j_4 < N_L - 1$, then the probability of drawing {L,L} is $C(N_L - j_3 - j_4,2)/C(N_S + N_L - k,2)$. The probability that L is observed and the unobserved L is assigned to B is 1-L. Then let

$$\alpha = (1\text{-}L)\, C(N_L - j_3 - j_4,2)/C(N_S + N_L - k,2).$$

If {S,L} has been drawn, then $j_1 + j_2 < N_S$ and $j_3 + j_4 < N_L$. If either $j_1 + j_2 = N_S$ or $j_3 + j_4 = N_L$ let $\beta = 0$. If not, the probability of {S,L} is $(N_S - j_1 - j_2)(N_L - j_3 - j_4)/(N_S + N_L - k,2) = \gamma$. The probability that S is observed and L is assigned to B is 0.5(1-S). Then let $\beta = 0.5(1\text{-}S)\gamma$.

Then $P(j_1,j_2+1,j_3,j_4) = (\alpha + \beta)\, P(j_1,j_2,j_3,j_4)$.

**Simulations**

The above recursion relations are difficult to solve – either by hand or on a small computer -- but the process is very easy to simulate. In the following simulation, S and L were randomly chosen between 0.05 and 0.90, and the initial number of coins of each species was randomly chosen between 1 and 20. 1000 trials were conducted for each. The column headed OK has 1 if $E(U\rightarrow A) < m < E(U\rightarrow B)$ and 0 otherwise. FOM = $(E(U\rightarrow B) - E(U\rightarrow A))/m$, and abbreviates Figure of Merit. The entire array was sorted smallest to largest on FOM. The results are presented in the table below, with a few observations following.

| S | #S | L | #L | m | A total | #A | E(U→A) | B total | #B | E(U→B) | OK | FOM |
|---|---|---|---|---|---|---|---|---|---|---|---|---|
| 0.03 | 17 | 0.22 | 21 | 0.135 | 332.24 | 2455 | 0.1353 | 2230.1 | 16545 | 0.1348 | 0 | -0.0037 |
| 0.26 | 18 | 0.35 | 6 | 0.2825 | 956.32 | 3377 | 0.2832 | 2435.3 | 8623 | 0.2824 | 0 | -0.0028 |
| 0.32 | 5 | 0.48 | 17 | 0.4436 | 2156.96 | 4856 | 0.4442 | 2722.4 | 6144 | 0.4431 | 0 | -0.0025 |
| 0.33 | 12 | 0.44 | 14 | 0.3892 | 1969 | 5053 | 0.3897 | 3090.01 | 7947 | 0.3888 | 0 | -0.0023 |
| 0.71 | 17 | 0.92 | 11 | 0.7925 | 8779.68 | 11067 | 0.7933 | 2321.83 | 2933 | 0.7916 | 0 | -0.0021 |
| 0.48 | 20 | 0.71 | 8 | 0.5457 | 4195.09 | 7688 | 0.5457 | 3439.39 | 6312 | 0.5449 | 0 | -0.0015 |
| 0.35 | 9 | 0.39 | 3 | 0.36 | 797.62 | 2214 | 0.3603 | 1362.38 | 3786 | 0.3598 | 0 | -0.0014 |
| 0.18 | 1 | 0.6 | 15 | 0.5737 | 2647.38 | 4623 | 0.5727 | 1931.7 | 3377 | 0.572 | 0 | -0.0012 |
| 0.36 | 3 | 0.48 | 5 | 0.435 | 771.36 | 1770 | 0.4358 | 971.04 | 2230 | 0.4354 | 0 | -0.0009 |
| 0.24 | 1 | 0.7 | 17 | 0.6744 | 4141.54 | 6134 | 0.6752 | 1933.52 | 2866 | 0.6746 | 0 | -0.0009 |
| 0.75 | 15 | 0.85 | 15 | 0.8 | 9554.65 | 11939 | 0.8003 | 2448.25 | 3061 | 0.7998 | 0 | -0.0006 |

| S | #S | L | #L | m | A total | #A | **E(U→A)** | B total | #B | **E(U→B)** | OK | FOM |
|---|---|---|---|---|---|---|---|---|---|---|---|---|
| 0.17 | 15 | 0.2 | 15 | 0.185 | 517.2 | 2793 | 0.1852 | 2259.06 | 12207 | 0.1851 | 0 | -0.0005 |
| 0.18 | 17 | 0.26 | 15 | 0.2175 | 739.06 | 3401 | 0.2173 | 2736.94 | 12599 | 0.2172 | 0 | -0.0005 |
| 0.8 | 18 | 0.83 | 14 | 0.8131 | 10550.78 | 12976 | 0.8131 | 2458.29 | 3024 | 0.8129 | 0 | -0.0002 |
| 0.41 | 13 | 0.44 | 11 | 0.4238 | 2159.94 | 5097 | 0.4238 | 2924.91 | 6903 | 0.4237 | 0 | -0.0002 |
| 0.09 | 11 | 0.11 | 9 | 0.099 | 96.68 | 978 | 0.0989 | 892.28 | 9022 | 0.0989 | 0 | 0.0000 |
| 0.42 | 13 | 0.42 | 21 | 0.42 | 3013.92 | 7176 | 0.42 | 4126.08 | 9824 | 0.42 | 1 | 0.0000 |
| 0.51 | 13 | 0.53 | 5 | 0.5156 | 2383.56 | 4624 | 0.5155 | 2256.24 | 4376 | 0.5156 | 1 | 0.0002 |
| 0.79 | 16 | 0.87 | 16 | 0.83 | 11013.12 | 13272 | 0.8298 | 2264.72 | 2728 | 0.8302 | 1 | 0.0005 |
| 0.4 | 6 | 0.58 | 12 | 0.52 | 2420.82 | 4671 | 0.5183 | 2245.14 | 4329 | 0.5186 | 0 | 0.0006 |
| 0.81 | 13 | 0.88 | 17 | 0.8497 | 10812.2 | 12726 | 0.8496 | 1933.22 | 2274 | 0.8501 | 1 | 0.0006 |
| 0.46 | 16 | 0.67 | 18 | 0.5712 | 5539.12 | 9727 | 0.5695 | 4147.57 | 7273 | 0.5703 | 0 | 0.0014 |
| 0.59 | 2 | 0.72 | 16 | 0.7056 | 4512.81 | 6394 | 0.7058 | 1841.87 | 2606 | 0.7068 | 0 | 0.0014 |
| 0.4 | 10 | 0.72 | 10 | 0.56 | 3144.56 | 5591 | 0.5624 | 2483.28 | 4409 | 0.5632 | 0 | 0.0014 |
| 0.8 | 10 | 0.86 | 18 | 0.8386 | 9835.66 | 11729 | 0.8386 | 1907.4 | 2271 | 0.8399 | 1 | 0.0016 |
| 0.75 | 2 | 0.84 | 6 | 0.8175 | 2639.07 | 3228 | 0.8176 | 632.19 | 772 | 0.8189 | 0 | 0.0016 |
| 0.79 | 6 | 0.9 | 14 | 0.867 | 7456.92 | 8603 | 0.8668 | 1212.86 | 1397 | 0.8682 | 1 | 0.0016 |
| 0.29 | 18 | 0.51 | 8 | 0.3577 | 1652.82 | 4632 | 0.3568 | 2991.68 | 8368 | 0.3575 | 0 | 0.0020 |
| 0.6 | 11 | 0.78 | 17 | 0.7093 | 7018.02 | 9894 | 0.7093 | 2918.1 | 4106 | 0.7107 | 1 | 0.0020 |
| 0.29 | 7 | 0.51 | 19 | 0.4508 | 2630.45 | 5835 | 0.4508 | 3236.59 | 7165 | 0.4517 | 1 | 0.0020 |
| 0.56 | 5 | 0.87 | 19 | 0.8054 | 7762.38 | 9631 | 0.806 | 1913.78 | 2369 | 0.8078 | 0 | 0.0022 |
| 0.14 | 14 | 0.19 | 18 | 0.1681 | 458.4 | 2730 | 0.1679 | 2233.65 | 13270 | 0.1683 | 1 | 0.0024 |
| 0.29 | 15 | 0.39 | 1 | 0.2963 | 690.65 | 2335 | 0.2958 | 1680.05 | 5665 | 0.2966 | 1 | 0.0027 |
| 0.28 | 5 | 0.37 | 5 | 0.325 | 536.09 | 1652 | 0.3245 | 1089.54 | 3348 | 0.3254 | 1 | 0.0028 |
| 0.75 | 8 | 0.84 | 20 | 0.8143 | 9315.06 | 11441 | 0.8142 | 2089.53 | 2559 | 0.8165 | 1 | 0.0028 |
| 0.72 | 10 | 0.87 | 4 | 0.7629 | 4101.66 | 5378 | 0.7627 | 1240.59 | 1622 | 0.7649 | 1 | 0.0029 |
| 0.25 | 1 | 0.76 | 19 | 0.7345 | 5356.85 | 7280 | 0.7358 | 2008.04 | 2720 | 0.7382 | 0 | 0.0033 |
| 0.77 | 10 | 0.89 | 10 | 0.83 | 6898.8 | 8316 | 0.8296 | 1402.04 | 1684 | 0.8326 | 1 | 0.0036 |
| 0.15 | 16 | 0.36 | 12 | 0.24 | 793.11 | 3319 | 0.239 | 2562.48 | 10681 | 0.2399 | 0 | 0.0038 |
| 0.32 | 11 | 0.68 | 17 | 0.5386 | 4047.76 | 7517 | 0.5385 | 3504.48 | 6483 | 0.5406 | 1 | 0.0039 |
| 0.1 | 1 | 0.44 | 21 | 0.4245 | 1949.18 | 4610 | 0.4228 | 2713 | 6390 | 0.4246 | 1 | 0.0042 |
| 0.48 | 12 | 0.61 | 8 | 0.532 | 2844.1 | 5357 | 0.5309 | 2476.16 | 4643 | 0.5333 | 1 | 0.0045 |
| 0.38 | 11 | 0.63 | 15 | 0.5242 | 3530.55 | 6735 | 0.5242 | 3298.95 | 6265 | 0.5266 | 1 | 0.0046 |
| 0.09 | 18 | 0.32 | 18 | 0.205 | 764.12 | 3742 | 0.2042 | 2926.11 | 14258 | 0.2052 | 1 | 0.0049 |
| 0.31 | 5 | 0.38 | 1 | 0.3217 | 305.99 | 955 | 0.3204 | 658.45 | 2045 | 0.322 | 1 | 0.0050 |
| 0.4 | 1 | 0.57 | 11 | 0.5558 | 1865.52 | 3365 | 0.5544 | 1468.97 | 2635 | 0.5575 | 1 | 0.0056 |
| 0.24 | 11 | 0.74 | 5 | 0.3963 | 1257.18 | 3182 | 0.3951 | 1914.82 | 4818 | 0.3974 | 1 | 0.0058 |
| 0.15 | 3 | 0.36 | 7 | 0.297 | 420.81 | 1425 | 0.2953 | 1062.09 | 3575 | 0.2971 | 1 | 0.0061 |
| 0.28 | 3 | 0.79 | 15 | 0.705 | 4443.71 | 6317 | 0.7035 | 1899.25 | 2683 | 0.7079 | 1 | 0.0062 |
| 0.36 | 8 | 0.65 | 14 | 0.5445 | 3259.99 | 6021 | 0.5414 | 2715.51 | 4979 | 0.5454 | 1 | 0.0073 |
| 0.34 | 16 | 0.77 | 4 | 0.426 | 1807.18 | 4263 | 0.4239 | 2451.1 | 5737 | 0.4272 | 1 | 0.0077 |
| 0.32 | 12 | 0.65 | 16 | 0.5086 | 3621.37 | 7166 | 0.5054 | 3482.46 | 6834 | 0.5096 | 1 | 0.0083 |
| 0.6 | 9 | 0.9 | 13 | 0.7773 | 6588.9 | 8517 | 0.7736 | 1937.1 | 2483 | 0.7801 | 1 | 0.0084 |
| 0.23 | 6 | 0.65 | 16 | 0.5355 | 3146.28 | 5904 | 0.5329 | 2739.52 | 5096 | 0.5376 | 1 | 0.0088 |
| 0.15 | 7 | 0.41 | 21 | 0.345 | 1676.65 | 4877 | 0.3438 | 3164.53 | 9123 | 0.3469 | 1 | 0.0090 |
| 0.08 | 5 | 0.48 | 19 | 0.3967 | 1872.64 | 4748 | 0.3944 | 2888.56 | 7252 | 0.3983 | 1 | 0.0098 |
| 0.47 | 3 | 0.76 | 13 | 0.7056 | 3972.85 | 5659 | 0.702 | 1659.68 | 2341 | 0.709 | 1 | 0.0099 |

| S | #S | L | #L | m | A total | #A | **E(U→A)** | B total | #B | **E(U→B)** | OK | FOM |
|---|---|---|---|---|---|---|---|---|---|---|---|---|
| 0.47 | 3 | 0.76 | 13 | 0.7056 | 3972.85 | 5659 | 0.702 | 1659.68 | 2341 | 0.709 | 1 | 0.0099 |
| 0.53 | 18 | 0.98 | 20 | 0.7668 | 11150.58 | 14586 | 0.7645 | 3409.52 | 4414 | 0.7724 | 1 | 0.0103 |
| 0.2 | 10 | 0.38 | 2 | 0.23 | 315.26 | 1381 | 0.2283 | 1065.82 | 4619 | 0.2307 | 1 | 0.0104 |
| 0.32 | 12 | 0.75 | 14 | 0.5515 | 3924.3 | 7133 | 0.5502 | 3262.47 | 5867 | 0.5561 | 1 | 0.0107 |
| 0.39 | 4 | 0.78 | 14 | 0.6933 | 4303.26 | 6221 | 0.6917 | 1943.37 | 2779 | 0.6993 | 1 | 0.0110 |
| 0.26 | 2 | 0.79 | 14 | 0.7238 | 4200.21 | 5838 | 0.7195 | 1572.83 | 2162 | 0.7275 | 1 | 0.0111 |
| 0.34 | 19 | 0.85 | 7 | 0.4773 | 2944.91 | 6200 | 0.475 | 3266.21 | 6800 | 0.4803 | 1 | 0.0111 |
| 0.35 | 17 | 0.98 | 15 | 0.6453 | 6641.46 | 10350 | 0.6417 | 3666.53 | 5650 | 0.6489 | 1 | 0.0112 |
| 0.49 | 20 | 0.94 | 2 | 0.5309 | 3102.75 | 5850 | 0.5304 | 2762.9 | 5150 | 0.5365 | 1 | 0.0115 |
| 0.48 | 8 | 0.76 | 6 | 0.6 | 2511.2 | 4198 | 0.5982 | 1695.52 | 2802 | 0.6051 | 1 | 0.0115 |
| 0.5 | 20 | 0.92 | 4 | 0.57 | 3861.48 | 6804 | 0.5675 | 2983.14 | 5196 | 0.5741 | 1 | 0.0116 |
| 0.14 | 3 | 0.37 | 7 | 0.301 | 438.64 | 1464 | 0.2996 | 1072.8 | 3536 | 0.3034 | 1 | 0.0126 |
| 0.28 | 15 | 0.7 | 9 | 0.4375 | 2240.56 | 5161 | 0.4341 | 3011.54 | 6839 | 0.4403 | 1 | 0.0142 |
| 0.22 | 7 | 0.57 | 11 | 0.4339 | 1641.06 | 3813 | 0.4304 | 2265.34 | 5187 | 0.4367 | 1 | 0.0145 |
| 0.6 | 5 | 0.94 | 13 | 0.8456 | 6391.34 | 7586 | 0.8425 | 1209.14 | 1414 | 0.8551 | 1 | 0.0149 |
| 0.1 | 20 | 0.72 | 18 | 0.3937 | 2922.38 | 7468 | 0.3913 | 4587.38 | 11532 | 0.3978 | 1 | 0.0165 |
| 0.11 | 7 | 0.53 | 19 | 0.4169 | 2201.77 | 5339 | 0.4124 | 3211.93 | 7661 | 0.4193 | 1 | 0.0166 |
| 0.35 | 10 | 0.65 | 6 | 0.4625 | 1716.3 | 3732 | 0.4599 | 1996.3 | 4268 | 0.4677 | 1 | 0.0169 |
| 0.14 | 18 | 0.65 | 10 | 0.3221 | 1434.17 | 4543 | 0.3157 | 3055.43 | 9457 | 0.3231 | 1 | 0.0230 |
| 0.04 | 10 | 0.25 | 16 | 0.1692 | 361.41 | 2184 | 0.1655 | 1837.12 | 10816 | 0.1699 | 1 | 0.0260 |
| 0.14 | 17 | 0.7 | 17 | 0.42 | 2992.92 | 7226 | 0.4142 | 4166.12 | 9774 | 0.4262 | 1 | 0.0286 |
| 0.2 | 15 | 0.74 | 7 | 0.3718 | 1489.7 | 4033 | 0.3694 | 2648.36 | 6967 | 0.3801 | 1 | 0.0288 |
| 0.2 | 14 | 0.55 | 8 | 0.3273 | 1142.25 | 3557 | 0.3211 | 2461.6 | 7443 | 0.3307 | 1 | 0.0293 |
| 0.24 | 17 | 0.84 | 3 | 0.33 | 1078.68 | 3332 | 0.3237 | 2226.12 | 6668 | 0.3339 | 1 | 0.0309 |
| 0.36 | 13 | 0.81 | 9 | 0.5441 | 3220.56 | 5996 | 0.5371 | 2772.54 | 5004 | 0.5541 | 1 | 0.0312 |
| 0.15 | 18 | 0.65 | 14 | 0.3687 | 2130.3 | 5932 | 0.3591 | 3734.2 | 10068 | 0.3709 | 1 | 0.0320 |
| 0.13 | 3 | 0.92 | 13 | 0.7719 | 4710.53 | 6148 | 0.7662 | 1465.26 | 1852 | 0.7912 | 1 | 0.0324 |
| 0.17 | 6 | 0.95 | 16 | 0.7373 | 5961.57 | 8181 | 0.7287 | 2128.15 | 2819 | 0.7549 | 1 | 0.0355 |
| 0.11 | 16 | 0.89 | 16 | 0.5 | 3917.37 | 7965 | 0.4918 | 4104.47 | 8035 | 0.5108 | 1 | 0.0380 |
| 0.05 | 13 | 0.16 | 3 | 0.0706 | 41.07 | 608 | 0.0675 | 519.42 | 7392 | 0.0703 | 0 | 0.0397 |
| 0.07 | 19 | 0.75 | 13 | 0.3463 | 1863.8 | 5536 | 0.3367 | 3694.56 | 10464 | 0.3531 | 1 | 0.0474 |
| 0.36 | 10 | 0.94 | 8 | 0.6178 | 3412.82 | 5644 | 0.6047 | 2130.94 | 3356 | 0.635 | 1 | 0.0490 |
| 0.07 | 15 | 0.49 | 9 | 0.2275 | 605.85 | 2739 | 0.2212 | 2152.71 | 9261 | 0.2324 | 1 | 0.0492 |
| 0.27 | 11 | 0.95 | 11 | 0.61 | 4033.36 | 6728 | 0.5995 | 2697.72 | 4272 | 0.6315 | 1 | 0.0525 |
| 0.02 | 15 | 0.37 | 11 | 0.1681 | 340.08 | 2129 | 0.1597 | 1833.37 | 10871 | 0.1686 | 1 | 0.0529 |
| 0.18 | 12 | 0.94 | 8 | 0.484 | 2275.32 | 4876 | 0.4666 | 2551 | 5124 | 0.4979 | 1 | 0.0647 |
| 0.01 | 14 | 0.95 | 12 | 0.4438 | 2452.23 | 5805 | 0.4224 | 3274.53 | 7195 | 0.4551 | 1 | 0.0737 |
| 0.07 | 11 | 0.53 | 7 | 0.2489 | 530.77 | 2299 | 0.2309 | 1674.73 | 6701 | 0.2499 | 1 | 0.0763 |
| 0.05 | 1 | 0.5 | 3 | 0.3875 | 277.85 | 760 | 0.3656 | 509.3 | 1240 | 0.4107 | 1 | 0.1164 |
| 0.05 | 18 | 0.72 | 4 | 0.1718 | 283.81 | 1817 | 0.1562 | 1631.65 | 9183 | 0.1777 | 1 | 0.1251 |
| 0.09 | 4 | 0.58 | 4 | 0.335 | 420.77 | 1365 | 0.3083 | 925.6 | 2635 | 0.3513 | 1 | 0.1284 |
| 0.05 | 8 | 0.57 | 4 | 0.2233 | 269.66 | 1358 | 0.1986 | 1080.22 | 4642 | 0.2327 | 1 | 0.1527 |
| 0.01 | 8 | 0.53 | 2 | 0.114 | 43.68 | 572 | 0.0764 | 539.32 | 4428 | 0.1218 | 1 | 0.3982 |
| 0.08 | 2 | 0.97 | 2 | 0.525 | 428.55 | 1096 | 0.391 | 596.53 | 904 | 0.6599 | 1 | 0.5122 |

If one looks at the last few entries – the ones with FOM > 10, one notices that both #S (the number of marbles of weight S initially in the bag) and #L are small in all but one case. This might be the result of the fact that the SAP in these cases results in selection probabilities that differ considerably from independence. Conversely, when #S and #L are large, the selection probabilities might mirror independent probabilities more closely.

There were a number of instances in which $E(U \to B) < E(U \to B)$, but these all had FOMs very close to 0, as is the case with independent selection (Cor. 3.2). The author suspects, but is unable to prove, that $E(U \to A) \leq E(U \to B)$ in all cases.